\documentclass[12pt,reqno]{amsart}

\usepackage{setspace}

\usepackage
[breaklinks=true,
hypertexnames=false]   
{hyperref}

\usepackage{breakurl}   

\usepackage{enumitem} 

\theoremstyle{definition}

\DeclareMathOperator{\N}{{\mathbb N}}

\DeclareMathOperator{\R}{{\mathbb R}}

\numberwithin{equation}{section}

\newcommand\astr{{{}^\ast\hspace{-2.5pt}\R}}

\newcommand{\hr} {{{}^{\mathfrak{h}}\hspace*{-2.3pt}\R}}

\newcommand{\st}{\textbf{st}}

\author[Bair]{Jacques Bair}\address{J. Bair, HEC-ULG, University of
Liege, 4000 Belgium}\email{j.bair@ulg.ac.be}

\author[Blaszczyk]{Piotr B\l{}aszczyk}\address{P. B\l{}aszczyk, Institute
of Mathematics, Pedagogical University of Cracow,
Poland}\email{pb@up.krakow.pl}

\author[Heinig]{Peter Heinig} \address{P. Heinig, Germany}
\email{peter.c.heinig@gmail.com}

\author[Kanovei]{Vladimir Kanovei} \address{V. Kanovei, IPPI RAS,
Moscow, Russia}\email{kanovei@googlemail.com}

\author[Katz]{Mikhail G. Katz}\address{M. Katz, Department of
Mathematics, Bar Ilan University, Ramat Gan 5290002
Israel}\email{katzmik@macs.biu.ac.il}

\author[McGaffey]{Thomas McGaffey}\address{T. McGaffey, Rice
University, US}\email{thomasmcgaffey@sbcglobal.net}


\begin{document}


\thispagestyle{empty}


\title[Cauchy's work on integral geometry]{Cauchy's work on integral
geometry, centers of curvature, and other applications of
infinitesimals}

\begin{abstract}
Like his colleagues de Prony, Petit, and Poisson at the \emph{Ecole
Polytechnique}, Cauchy used infinitesimals in the Leibniz--Euler
tradition both in his research and teaching.  Cauchy applied
infinitesimals in an 1826 work in differential geometry where
infinitesimals are used neither as variable quantities nor as
sequences but rather as \emph{numbers}.  He also applied
infinitesimals in an 1832 article on integral geometry, similarly as
numbers.  We explore these and other applications of Cauchy's
infinitesimals as used in his textbooks and research articles.

An attentive reading of Cauchy's work challenges received views on
Cauchy's role in the history of analysis and geometry.  We demonstrate
the viability of Cauchy's infinitesimal techniques in fields as
diverse as geometric probability, differential geometry, elasticity,
Dirac delta functions, continuity and convergence.

\end{abstract}

\maketitle

\tableofcontents

\section{Introduction}

Cauchy was one of the founders of rigorous analysis.  However, the
meaning of \emph{rigor} to Cauchy is subject to debate among scholars.
Cauchy used the term \emph{infiniment petit} (infinitely small) both
as an adjective and as a noun, but the meaning of Cauchy's term is
similarly subject to debate.  While Judith Grabiner and some other
historians feel that a Cauchyan infinitesimal is a sequence tending to
zero (see e.g., \cite{Gr81}, 1981), others argue that there is a
difference between null sequences and infinitesimals in Cauchy; see
e.g., Laugwitz (\cite{La87}, 1987), Katz--Katz (\cite{11b}, 2011),
Borovik--Katz (\cite{12b}, 2012), Smory\'nski (\cite{Sm12}, 2012,
pp.\;361--373 and \cite{Sm17}, 2017, pp.\;56, 61), Bair et
al.\;(\cite{17a}, 2017 and \cite{19a}, 2019).

Cauchy used infinitesimals in the Leibniz--Euler tradition both in his
research and teaching, like his colleagues de Prony, Petit, and
Poisson (see Section~\ref{s64}).  In the present text we will examine
several applications Cauchy makes of infinitesimals, and argue that he
uses them as atomic entities (i.e., entities not analyzable into
simpler constituents) rather than sequences.  We explore Cauchy's use
of infinitesimals in areas ranging from Dirac delta to integral
geometry.

\section{Dirac delta, summation of series}
\label{s43}

We consider Cauchy's treatment of (what will be called later) a Dirac
delta function.  Cauchy explicitly uses a unit-impulse, infinitely
tall, infinitely narrow delta function, as an integral kernel.  Thus,
in 1827, Cauchy used infinitesimals in his definition of a Dirac delta
function:
\begin{quote}
Moreover one finds, denoting by~$\alpha$,~$\epsilon$ two infinitely
small numbers,%
\footnote{As discussed in (Laugwitz \cite{La89}, 1989), a further
condition needs to be imposed on~$\alpha$ and~$\epsilon$ in modern
mathematics to ensure the correctness of the formula.}
\begin{equation}
\label{151}
\frac{1}{2} \int_{a-\epsilon}^{a+\epsilon} F(\mu) \frac{\alpha \;
d\mu}{\alpha^2 + (\mu-a)^2} = \frac{\pi}{2} F(a)
\end{equation}
(Cauchy \cite{Ca27}, 1827, p.\;289; counter\;\eqref{151} added)
\end{quote}
A formula equivalent to~\eqref{151} was proposed by Dirac a century
later.%
\footnote{The key property of the Dirac delta ``function''~$\delta(x)$
is exemplified by the defining formula~$\int_{-\infty}^{\infty}
f(x)\delta(x)=f(0)$, where~$f(x)$ is any continuous function of~$x$
(Dirac \cite{Dir}, 1930/1958, p.\;59).}
The expression
\begin{equation}
\label{31b}
\frac{\alpha}{\alpha^2 + (\mu-a)^2}
\end{equation}
occurring in Cauchy's formula is known as the {\em Cauchy
distribution\/} in probability theory.  Here Cauchy specifies a
function which meets the criteria as set forth by Dirac a century
later.  Cauchy integrates the function~$F$ against the
kernel~\eqref{31b} as in formula~\eqref{151} so as to extract the
value of~$F$ at the point~$a$, exploiting the characteristic property
of a delta function.

From a modern viewpoint, formula~\eqref{151} holds up to an
infinitesimal error.  For obvious reasons, Cauchy was unfamiliar with
modern set-theoretic foundational \emph{ontology} of analysis (with or
without infinitesimals), but his \emph{procedures} find better proxies
in modern infinitesimal frameworks than Weierstrassian ones.%
\footnote{On the procedures/ontology distinction see B\l aszczyk et
al.\;(\cite{17d}, 2017).}
From the modern viewpoint, the right hand side of \eqref{151}, which
does not contain infinitesimals, is the \emph{standard part} (see
Section~\ref{s71}) of the left hand side, which does contain
infinitesimals.  Thus, a Cauchy distribution with an infinitesimal
scale parameter~$\alpha$ produces an entity with Dirac-delta behavior,
exploited by Cauchy already in 1827; see Katz--Tall (\cite{13g}, 2013)
for details.

Similarly, in his article (Cauchy \cite{Ca53}, 1853) on the
convergence of series of functions, infinitesimals are handled as
atomic inputs to functions.  Here Cauchy studies the series
\[
u_0+\ldots+ u_n +\ldots
\]
Cauchy proceeds to choose ``une valeur infiniment grande'' (an
infinitely large value) for the index~$n$ in \cite[p.\;456]{Ca53}.  He
then states his convergence theorem modulo a hypothesis that the
sum~$u_n+u_{n+1}+\ldots+u_{n'-1}$ should be
\begin{quote}
toujours infiniment petite pour des valeurs infiniment grandes des
\emph{nombres} entiers~$n$ et~$n'>n$ {\dots}%
\footnote{Translation: ``always infinitely small for infinitely large
values of whole numbers~$n$ and~$n'>n$.''}
(Cauchy \cite{Ca53}, 1853, p.\;457; emphasis added).
\end{quote}
Cauchy's proof of the continuity of the sum exploits the condition
that the sum~$u_n+u_{n+1}+\ldots+u_{n'-1}$ should be infinitesimal for
atomic infinitesimal inputs.  Cauchy writes down such an input in the
form~$x=\frac{1}{n}$; see e.g., \cite[p.\;457]{Ca53}.  For further
details on \cite{Ca53} see Section~\ref{s63b}; see also Bascelli et
al.\;(\cite{18e}, 2018).

\section{Differentials, infinitesimals, and derivatives}
\label{s21}

In his work in analysis, Cauchy carefully distinguishes between
differentials~$ds,dt$ which to Cauchy are noninfinitesimal variables,
on the one hand, and infinitesimal increments~$\Delta s, \Delta t$, on
the other:
\begin{quote}
[S]oit~$s$ une variable distincte de la variable primitive~$t$.  En
vertu des d\'efinitions adopt\'ees, le rapport entre les
diff\'erentielles~$ds$,~$dt$, sera la limite du rapport entre les
\emph{accroissements infiniment petits}~$\Delta s$,~$\Delta t$.%
\footnote{Let~$s$ be a variable distinct from the primitive variable
$t$.  By virtue of the definitions chosen, the ratio between the
differentials~$ds$,~$dt$ will be the limit of the ratio between the
infinitely small increments~$\Delta s$,~$\Delta t$.''}
(Cauchy \cite{Ca44}, 1844, p.\;11; emphasis added)
\end{quote}
Cauchy goes on to express such a relation by means of a formula in
terms of the infinitesimals~$\Delta s$ and~$\Delta t$:
\begin{quote}
On aura donc%
\footnote{Translation: ``One will then have.''}
\begin{equation}
\label{e31b}
\frac{ds}{dt} = \; \text{lim.}\, \frac{\Delta s}{\Delta t}
\end{equation}
(ibid., equation (1); significantly, the period after lim as in
``lim.'' is in the original; counter \eqref{e31b} added)
\end{quote}
In modern infinitesimal frameworks, the passage from the ratio of
infinitesimals such as~$\frac{\Delta s}{\Delta t}$ to the value of the
derivative is carried out by the \emph{standard part function}; see
equations~\eqref{e61} and~\eqref{e62} in Section~\ref{s71}.
Paraphrazing Cauchy's definition of the derivative as in \eqref{e31b}
in Archimedean terms would necessarily involve elements that are
inexplicit in the original definition.  Thus Cauchy's ``lim.''\;finds
a closer proxy in the notion of standard part, as in formula
\eqref{e83}, than in any notion of limit in the context of an
Archimedean continuum; see also Bascelli et al.\;(\cite{14a}, 2014).

\section{Integral geometry}
\label{s3}

An illuminating use of infinitesimals occurs in Cauchy's article in a
field today called \emph{integral geometry} (also known as
\emph{geometric probability}); see Hyk\v sov\'a et al.\;(\cite{Hy12},
2012, pp.\;3--4) for a discussion.

\subsection{Decomposition into infinitesimal segments}

Cauchy proved a formula known today as the Cauchy--Crofton formula%
\footnote{Or the Crofton formula; see e.g., Tabachnikov (\cite{Ta05},
2005, p.\;37).}
in his article (\cite{Ca50}, 1850; originally presented as
\cite{Ca32}, 1832).  Here Cauchy exploits a decomposition of a curve
into infinitesimal length elements (respectively, of a surface into
infinitesimal area elements) in an essential way in proving a formula
for the length of a plane curve (respectively, area of a surface
in~$3$-space).  Thus, Cauchy proves the formula
\begin{equation}
\label{e41b}
S= \frac{1}{4} \int_{-\pi}^{\pi} A\, dp
\end{equation}
for the length of a plane curve, in his Th\'eor\`eme I in
\cite[p.\;167--168]{Ca50}.  In formula~\eqref{e41b},~$p$ is the polar
angle (usually denoted~$\theta$ today), whereas~$A$ is the sum of the
orthogonal projections of the length elements onto a rotating line
with parameter~$p$.  Note that this is an exact formula (rather than
an approximation), typical of modern \emph{integral geometry}.

In his Th\'eor\`eme II, Cauchy goes on to prove a constructive
version, or a discretisation, of his Th\'eor\`eme I.\, Here Cauchy
replaces \emph{integrating} with respect to the variable-line
differential $dp$, by \emph{averaging} over a system of~$n$ equally
spaced lines (i.e., such that successive lines form equal angles).
Cauchy then obtains the approximation
\begin{equation}
\label{e41}
S= \frac{\pi}{2} M
\end{equation}
where~$M$ is the average.  Here the equality sign appears in
\cite[p.\;169]{Ca50} as in our formula~\eqref{e41}, and denotes
approximation.  Cauchy also provides an explicit error bound of
\begin{equation}
\label{e31}
\frac{\pi M}{2n^2}
\end{equation}
for the approximation, in \cite[p.\;169]{Ca50}.  Cauchy first proves
the result for a straight line segment, and then writes:
\begin{quote}
Le th\'eor\`eme II \'etant ainsi d\'emontr\'e pour le cas particulier
o\`u la quantit\'e~$S$ se r\'eduit \`a une longueur rectiligne~$s$, il
suffira, pour le d\'emontrer dans le cas contraire, de
d\'ecomposer~$S$ en \emph{\'el\'ements infiniment petits}.%
\footnote{Translation: ``Theorem II having been proved for the special
case when the quantity~$S$ is a straight line segment~$s$, it would be
sufficient, to prove it in the contrary [i.e., the general] case, to
decompose~$S$ into infinitely small elements.''}
(Cauchy \cite{Ca50}, 1850, p.\;171; emphasis added)
\end{quote} 
Thus Cauchy obtains a sequence of error bounds of the form~\eqref{e31}
that improve (become smaller) as~$n$ increases.

\subsection{Analysis of Cauchy's argument}

Cauchy exploits two entities which need to be carefully distinguished
to keep track of the argument:
\begin{enumerate}
\item
the curve itself, and
\item
the circle (or in modern terminology, the Grassmannian) of directions
parametrized by~$p$ (counterpart of the modern polar angle~$\theta$).
\end{enumerate}
Note that Cauchy treats the curve and the Grassmannian differently.
Namely, the curve is subdivided into infinitely many infinitesimal
elements of length.  Meanwhile, as far as the Grassmannian is
concerned, Cauchy works with a finite~$n$, chooses~$n$ directions that
are equally spaced, and is interested in the asymptotic behavior of
the \emph{sequence} of error estimates~\eqref{e31} as~$n$ tends to
infinity.

If an infinitesimal merely meant a variable quantity or sequence to
Cauchy, then there shouldn't be any difference in Cauchy's treatment
of the curve and the Grassmannian; both should be sequences.  However,
Cauchy does treat them differently:
\begin{itemize}
\item
the curve is viewed as an aggregate of infinitely many infinitesimal
elements;
\item
the circle of directions is decomposed into~$n$ segments and the focus
is on the asymptotic behavior of the error bound as a function of~$n$.
\end{itemize}
The approach in the 1832/1850 paper on integral geometry indicates
that Cauchy's infinitesimal (the element of length decomposing the
curve) is not a variable quantity or a sequence, but rather an atomic
entity, as discussed in Section~\ref{s43}.

\section{Centers of curvature, elasticity}

In studying the geometry of curves, Cauchy routinely exploits
infinitesimals and related notions such as infinite proximity.  We
will analyze Cauchy's book \emph{Le\c cons sur les applications du
calcul infinit\'esimal \`a la g\'eom\'etrie} (\cite{Ca26}, 1826).

\subsection{\emph{Angle de contingence}}

Cauchy starts by defining the \emph{angle de
contingence}~$\pm\Delta\tau$ as the angle between the two tangent
lines of an arc~$\pm\Delta s$ at its extremities.%
\footnote{To follow the mathematics it is helpful to think of
parameter~$\tau$ as the angle measured counterclockwise between the
positive direction of the~$x$-axis and the tangent vector to the
curve.}
He then considers the normals to the curve at the extremities of the
arc starting at the point~$(x,y)$.

\subsection{Center of curvature and radius of curvature}

Cauchy goes on to give \emph{two} definitions of both the center of
curvature and the radius of curvature.  Thus, he writes:
\begin{quote}
la distance du point~$(x, y)$ au point de rencontre des deux normales
est sensiblement \'equivalente au rayon d'un cercle qui aurait la
m\^eme courbure que la courbe.%
\footnote{Translation: ``the distance from the point $(x,y)$ to the
intersection point of two normals is appreciably equivalent to the
radius of a circle which would have the same curvature as the
curve.''}
(Cauchy \cite{Ca26}, 1826, p.\;98)
\end{quote}
Notice that Cauchy mentions two items:
\begin{enumerate}
[label={(Ca\theenumi)}]
\item
\label{i1b}
the intersection of the two normal lines produces a point which will
generate the \emph{center} of curvature, and
\item
the distance between~$(x,y)$ and the point defined in item~\ref{i1b}
which will generate the \emph{radius} of curvature.
\end{enumerate}
Thus the \emph{radius} of curvature is naturally defined in terms of
the \emph{center} of curvature (namely, as the distance between the
point~$(x,y)$ and the center of curvature).

\subsection{$\varepsilon$, \emph{nombre infiniment petit}}
\label{s63}

Next, Cauchy chooses an infinitesimal number~$\varepsilon$ and
exploits the law of sines to write down a relation that will give an
expression for the radius of curvature:
\begin{quote}
[S]i l'on d\'esigne par~$\varepsilon$ un \emph{nombre} infiniment
petit, on aura%
\footnote{Translation: ``If one denotes by~$\varepsilon$ an infinitely
small number, one will obtain''}
\begin{equation}
\label{e61b}
\frac{\sin\left(\frac{\pi}{2}\pm\varepsilon\right)}{r}=
\frac{\sin(\pm\Delta\tau)}{\sqrt{\Delta x^2+\Delta y^2}}
\end{equation}
(Cauchy \cite{Ca26}, 1826, p.\;98; emphasis and counter
``\eqref{e61b}'' added).
\end{quote}
Note that Cauchy describes his infinitesimal~$\varepsilon$ neither as
a \emph{sequence} nor as a \emph{variable quantity} but rather as an
infinitely small \emph{number} (``nombre'').  At the next stage,
Cauchy passes to the limit to obtain:
\begin{quote}
On en conclura, en passant aux \emph{limites},
\begin{equation}
\label{e62b}
\frac{1}{\rho} = \pm \frac{d\tau}{\sqrt{dx^2+dy^2}}
\end{equation}
\cite[p.\;99]{Ca26} (emphasis and counter ``\eqref{e62b}'' added)
\end{quote}
It is instructive to analyze what happens exactly in passing from
formula~\eqref{e61b} to formula~\eqref{e62b}.  Here Cauchy replaces
the infinitesimals~$\Delta x$, $\Delta y$, and~$\sin\Delta\tau$ by the
corresponding differentials~$dx$,~$dy$, and~$d\tau$.  The
expression~$\sin(\frac{\pi}{2}\pm\varepsilon)$ is infinitely close
to~$1$ whereas~$r$ is infinitely close to~$\rho$, justifying the
replacements in the left-hand side of Cauchy's equation.

As in Cauchy's definition of derivative analyzed in Section~\ref{s21},
Cauchy's \emph{limite} here admits of a close proxy in the standard
part function (see Section~\ref{s71}).  Meanwhile, any attempt to
interpret Cauchy's procedure in the context of an Archimedean
continuum will have to deal with the nettlesome issue of the absence
of Cauchyan infinitesimals like~$\Delta x$,~$\Delta y$, 
$\Delta\tau$, and~$\varepsilon$ in such a framework.

\subsection
{Second characterisation of radius and center of curvature}

Using his formula for~$\rho$, Cauchy goes on to give his \emph{second}
characterisation of the radius of curvature and center of curvature:
\begin{quote}
Ce rayon, port\'e \`a partir du point~$(x, y)$ sur la normale qui
renferme ce point, est ce qu'on nomme le \emph{rayon de courbure} de
la courbe propos\'ee, relatif au point dont il s'agit, et l'on appelle
\emph{centre de courbure} celle des extr\'emit\'es du rayon de
courbure que l'on peut consid\'erer comme le point de rencontre de
deux normales infiniment voisines.  (Cauchy \cite{Ca26}, p.\;99,
emphasis in the original)
\end{quote}
Here Cauchy notes that the center of curvature can be viewed as the
other endpoint of the vector of length~$\rho$ starting at the
point~$(x,y)$ and normal to the curve.  Cauchy then reiterates the
earlier definition of the center of curvature of a plane curve in
terms of the intersection point of a pair of \emph{infinitely close}
normals, as Leibniz may have done; see Katz--Sherry (\cite{13f},
2013).  Note that neither the center of curvature nor the radius of
curvature are defined using a notion of limit in the context of an
Archimedean continuum.  

Cauchy's presentation of infinitesimal techniques here contains no
trace of the variable quantities or sequences exploited in his
textbooks in the definitions of infinitesimals.  To adapt Cauchy's
definition of center of curvature to modern custom, it is certainly
possible to paraphrase it in the context of an Archimedean continuum.
This can be done for example by taking a suitable sequence of (pairs
of) normals and passing to a limit.  However, such a paraphrase would
not be faithful to Cauchy's own \emph{procedure}.  We will follow the
continuation of Cauchy's analysis in Section~\ref{s55}.

\subsection
{Formula for~$\rho$ in terms of infinitesimal displacements}
\label{s55}

At this stage, Cauchy's goal is to develop a formula for the radius of
curvature~$\rho$ of the curve at a point~$(x,y)$.  Cauchy seeks to
express~$\rho$ in terms of the distance between the pair of points
obtained from~$(x,y)$ by means of equal infinitesimal displacements,
one along the curve and the other along the tangent.  To this end,
Cauchy starts by choosing an infinitesimal~$i$:
\begin{quote}
Ajoutons que, si, \`a partir du point~$(x,y)$, on porte sur la courbe
donn\'ee et sur sa tangente, prolong\'ees dans le m\^eme sens que
l'arc~$s$, des longueurs \'egales et infiniment petites
repr\'esent\'ees par~$i$, on trouvera, pour les coordonn\'ees de
l'extr\'emit\'e de la seconde longueur,
\[
x+ i\frac{dx}{ds}, \quad y+i\frac{dy}{ds}
\]
et, pour les coordonn\'ees de l'extr\'emit\'e de la premi\`ere,
\[
x+i\frac{dx}{ds} +\frac{i^2}{2}\left(\frac{d^2x}{ds^2}+I\right), \quad
y+i\frac{dy}{ds} +\frac{i^2}{2}\left(\frac{d^2y}{ds^2}+J\right), 
\]
$I,J$ d\'esignant des quantit\'es infiniment petites.  (Cauchy
\cite{Ca26}, 1826, p.\;105)
\end{quote}
Cauchy refers to the endpoints of the infinitesimal segment of the
curve itself as the \emph{extremities}.  He denotes the distance
between the two extremities by~$\gamma$.%
\footnote{Actually Cauchy uses a slightly different symbol not
available in modern fonts.}
Then straightforward calculations produce the following formula
for~$\rho$:
\begin{quote}
De cette derni\`ere formule {\ldots} on tire
\begin{equation}
\label{e63}
\rho=\lim \frac{i^2}{2\gamma}.
\end{equation}
(ibid.; counter ``\eqref{e63}'' added)
\end{quote}
Note that in Cauchy's formula that we labeled~\eqref{e63}, the
symbol~``$\lim$'' is applied to a ratio of two infinitesimals.
Therefore the use of~$\lim$ here is analogous to the use of the
standard part function as in~\eqref{e83}.  Cauchy employed a similar
technique in the definition of derivative analyzed in
Section~\ref{s21}, and in passing from formula \eqref{e61b} to
formula~\eqref{e62b} as analyzed in Section~\ref{s63}.  Cauchy
concludes as follows:
\begin{quote}
En cons\'equence, \emph{pour obtenir le rayon de courbure d'une courbe
en un point donn\'e, il suffit de porter sur cette courbe et sur sa
tangente, prolong\'ees dans le m\^eme sens, des longueurs \'egales et
infiniment petites, et de diviser le carr\'e de l'une d'elles par le
double de la distance comprise entre les deux extr\'emit\'es}.  La
limite du quotient est la valeur exacte du rayon de courbure.
(op.\;cit., pp.\;105--106; emphasis in the original)
\end{quote}
Cauchy's \emph{limite} here again plays the role of the standard part
\eqref{e61}.

\subsection{Elasticity}

Another example of Cauchy's application of infinitesimals is his
foundational article on elasticity (\cite{Ca23b}, 1823) where \emph{un
\'el\'ement infiniment petit} is exploited on page 302.  The article
is mentioned by Freudenthal in (\cite{Fr71}, 1971, p.\;378); for
details see Belhoste (\cite{Be91}, 1991, p.\;94).

\section{Continuity in Cauchy}
\label{six}

In his \emph{Cours d'Analyse} (\cite{Ca21}, 1821), Cauchy comments as
follows concerning the continuity of a few functions, including the
function~$\frac{a}{x}$ in the range between~$0$ and infinity:
\begin{quote}
[E]ach of these functions is continuous in the neighborhood of any
finite value given to the variable~$x$ if that finite value is
contained \ldots, for the function~$\frac{a}{x}$ \ldots{}, between the
limits~$x=0$ and~$x=\infty$.  (Cauchy as translated by Bradley and
Sandifer%
\footnote{Reinhard Siegmund-Schultze writes: ``By and large, with few
exceptions to be noted below, the translation is fine'' (\cite{Si},
2009).}
\cite{BS}, 2009, p.\;27)
\end{quote}
Here Cauchy asserts the continuity of the function~$\frac{a}{x}$ for
finite values of~$x$ contained between the limits~$0$ and~$\infty$.

\subsection{Ambiguity in definition of continuity}
\label{s61}

From a modern viewpoint, the above definition is ambiguous.
Interpreting it as continuity on~$(0,\infty)$ would rule out the
possibility of interpreting Cauchy's continuity as \emph{uniform
continuity}, since~$\frac{1}{x}$ is not uniformly continuous
on~$(0,\infty)$.  Interpreting it as saying that for every real~$x$
there is a neighborhood of~$x$ where the function is continuous, would
not rule out a \emph{uniform} interpretation; e.g., the
function~$\frac1x$ is uniformly continuous in a suitable neighborhood
of each nonzero real point.

Cauchy on occasion mentions that~$x$ is a \emph{real} (as opposed to
\emph{complex}) variable.  However, identifying Cauchy's notion with
the modern notion of real number would clearly be problematic.  Cauchy
seems not to have elaborated a distinction between finer types of
continuity we are familiar in modern mathematics, such as ordinary
pointwise continuity \emph{vs} uniform continuity.

%
%

\subsection{From variables to infinitesimals}
\label{s52}

There was a transformation in Cauchy's thinking about continuity from
an 1817 treatment in terms of \emph{variables} to an 1821 treatment in
terms of \emph{infinitesimals}.  In 1817, Cauchy defined continuity of
$f$ in terms of commutation of taking limit and evaluating the
function:
\begin{quote}
\emph{La limite d'une fonction continue de plusieurs variables est la
m\^eme fonction de leur limite}.  Cons\'equence de ce Th\'eor\`eme
relativement \`a la continuit\'e des fonctions compos\'ees qui ne
d\'ependent que d'une seule variable.%
\footnote{Translation: ``The limit of a continuous function of several
variables is [equal to] the same function of their limit.
Consequences of this Theorem with regard to the continuity of
composite functions dependent on a single variable.''  The reference
for this particular lesson in the Archives of the Ecole Polytechnique
is as follows: Le 4 Mars 1817, la le\c con 20.  Archives E. P., X II
C7, Registre d'instruction 1816--1817.}
(Cauchy as quoted by Guitard
\cite{Gu86}, 1986, p.\;34; emphasis added; cf.\;Belhoste \cite{Be91},
1991, p.\;309)
\end{quote}
Four years later in his \emph{Cours d'Analyse}, Cauchy defined
continuity as follows:
\begin{quote}
Among the objects related to the study of infinitely small quantities,
we ought to include ideas about the continuity and the discontinuity
of functions.  In view of this, let us first consider functions of a
single variable.  Let~$f(x)$ be a function of the variable~$x$, and
suppose that for each value of~$x$ between two given limits, the
function always takes a unique finite value.  If, beginning with a
value of~$x$ contained between these limits, we add to the
variable~$x$ an \emph{infinitely small increment}~$\alpha$, the
function itself is incremented by the difference~$f(x+\alpha)-f(x)$,
which depends both on the new variable~$\alpha$ and on the value
of~$x$.  Given this, the function~$f(x)$ is a \emph{continuous}
function of~$x$ between the assigned limits if, for each value of~$x$
between these limits, the numerical value of the
difference~$f(x+\alpha)-f(x)$ decreases indefinitely with the
numerical value of~$\alpha$.  (Cauchy as translated in
\cite[p.\;26]{BS}; emphasis on ``continuous'' in the original;
emphasis on ``infinitely small increment'' added)
\end{quote}
This definition can be thought of as an intermediary one between the
1817 definition purely in terms of variables, and his second 1821
definition stated purely in terms of infinitesimals.  Cauchy's second
definition summarizes the definition just given as follows:
\begin{quote}
In other words, \emph{the function~$f(x)$ is continuous with respect
to~$x$ between the given limits if, between these limits, an
infinitely small increment in the variable always produces an
infinitely small increment in the function itself}.  (ibid.; emphasis
in the original)
\end{quote}
Cauchy's second definition just quoted can be compared with one of the
first modern ones; see formula~\eqref{e74}.  Cauchy concludes his
discussion of continuity and discontinuity as follows:
\begin{quote}
We also say that the function~$f(x)$ is a continuous function of the
variable~$x$ in a neighborhood of a particular value of the
variable~$x$ whenever it is continuous between two limits of~$x$ that
enclose that particular value, even if they are very close together.
Finally, whenever the function~$f(x)$ ceases to be continuous in the
neighborhood of a particular value of~$x$, we say that it becomes
discontinuous, and that there is \emph{solution%
\footnote{meaning \emph{dissolution}, i.e., absence (of continuity).}
of continuity} for this particular value. (ibid.; emphasis in the
original)
\end{quote}
Three salient points emerge from these passages:
\begin{enumerate}
\item
Cauchy makes it clear at the outset that in his mind continuity is
``among the objects related to the study of infinitely small
quantities'';
\item
the infinitely small~$\alpha$ is used conspicuously in the
definitions;
\item
conspicuously absent from Cauchy's multiple definitions of continuity
is the notion of \emph{limit}.%
\footnote{The word \emph{limit} itself does occur in Cauchy's
definitions here but in an entirely different sense of \emph{endpoint
of an interval} where inputs to the function originate (what we would
call today the \emph{domain} of the function); cf.\;Smory\'nski
(\cite{Sm17}, 2017, p.\;52, note 48).}
\end{enumerate}
Yushkevich observes in this connection that ``the definition of
continuity in Cauchy is as far from the \emph{Epsilontik} as his
definition of limit'' (\cite{Yu86}, 1986, p.\;69).

\subsection{The 1853 definition}
\label{s63b}

Some three decades later in 1853, Cauchy defined continuity in a
similar fashion:
\begin{quote}
\ldots une fonction~$u$ de la variable r\'eelle~$x$ sera
\emph{continue}, entre deux limites donn\'ees de~$x$, si, cette
fonction admettant pour chaque valeur interm\'ediaire de~$x$ une
valeur unique et finie, un accroissement infiniment petit attribu\'e
\`a la variable produit toujours, entre les limites dont il s'agit, un
accroissement infiniment petit de la fonction elle-m\^eme.  (Cauchy
\cite{Ca53}, 1853; emphasis in the original)
\end{quote}
Cauchy's 1853 definition echoes the 1821 definition given in
Section~\ref{s52}, where Cauchy denoted his infinitely small~$\alpha$
and required the difference~$f(x+\alpha)-f(x)$ to be infinitesimal as
a criterion for continuity of the function~$f$.

Cauchy's 1821 example of the function~$\frac{1}{x}$ between~$0$ and
infinity suggests that Cauchy's definition of continuity is, from a
modern viewpoint, somewhat ambiguous, as discussed in
Section~\ref{s61}.  Resolving the ambiguity by attributing uniform
continuity to Cauchy may not preserve such inherent ambiguity.

A possible interpretation of Cauchy's comments is available in the
context of an infinitesimal-enriched continuum.  Here one can
interpret~$x$ as referring to an assignable value (i.e., what we refer
to today as a \emph{real} value), and~$\alpha$ an (inassignable)
infinitesimal.  Then~$\frac{1}{x}$ is continuous in a neighorbood
of~$x$ in the sense that for each infinitesimal~$\alpha$ the
difference~$f(x+\alpha)-f(x)$ is also infinitesimal.

\subsection{Contingency and determinacy}

We wish to suggest, following Hacking (\cite{Ha14}, 2014,
pp.\;72--75), the possibility of alternative courses for the
development of analysis (a \emph{Latin model} as opposed to a
\emph{butterfly model}).%
\footnote{Hacking contrasts a model of a deterministic biological
development of animals like butterflies (the
egg--larva--cocoon--butterfly sequence), as opposed to a model of a
contingent historical evolution of languages like Latin.}
From such a standpoint, the traditional assumption that the historical
development inexorably led to modern classical analysis (as formalized
by Weierstrass and others) remains merely a hypothesis.  A reader of
Dani--Papadopoulous \cite{Pa19} may be surprised to learn that
\begin{quote}
Cauchy gave a faultless definition of continuous function, using the
notion of `limit' for the first time.  Following Cauchy's idea,
Weierstrass popularized the $\epsilon$-$\delta$ argument in the
1870's.  (\cite{Pa19}, 2019, p.\;283)
\end{quote}
Such views fit well with a deterministic \emph{butterfly} model
leading from Cauchy to Weierstrass.  However, such views are not
merely anachronistic but contrary to fact, as we saw in
Sections~\ref{s52} and \ref{s63b}.  Cauchy did write: ``Lorsque les
valeurs num\'eriques successives d'une m\^eme variable d\'ecroissent
ind\'efiniment, de mani\`ere \`a s'abaisser au-dessous de tout nombre
donn\'e, cette variable devient ce qu'on nomme un \emph{infiniment
petit} ou une quantit\'e \emph{infiniment petite}.  Une variable de
cette esp\`ece a z\'ero pour limite'' \cite[p.\;4]{Ca21} (emphasis in
the original).  However, interpreting Cauchy's wording as an
anticipation of the modern \emph{Epsilontik} notion of limit would be
anachronistic, since Cauchy's wording here echoes formulations
provided by his teacher Lacroix, and even earlier formulations found
in Leibniz, not to speak of the ancient method of exhaustion; see Bair
et al.\;(\cite{19a}, 2019) for details.

We argue, following Robinson (\cite{Ro66}, 1966) and Laugwitz
\cite{La87}, 1987), that the procedures of Leibniz, Euler and Cauchy
were closer to the procedures in Robinson's framework than the
procedures in a Weierstrassian framework.  On this view,
interpretation of the work of Leibniz, Euler, and Cauchy in analysis
is more successful in a modern infinitesimal framework than a modern
Archimedean one; see e.g., Bair et al.\;(\cite{18a}, 2018) on Leibniz
and Bair et al.\;(\cite{17b}, 2017) on Euler.  For a survey of
infinitesimal mathematics and its history see e.g., Robinson
(\cite{Ro66}, 1966, chapter\;10, pp.\;260--282).

\section
{Reception of Cauchy's ideas among his colleagues}
\label{s64}

For historians advocating an externalist approach to the history of
mathematics, it is important to consider the reception of Cauchy's
ideas among his contemporaries.  Cauchy's contemporaries and
colleagues at the \emph{Ecole Polytechnique} (Poisson and others; see
below) had specific ideas about what \emph{infinitely small} meant.
One cannot provide a proper analysis of Cauchy's notion without taking
into account the ideas on the subject among his contemporaries.  There
seems to be little reason to doubt that the notion of infinitely small
in the minds of Poisson, de Prony, Petit, and others was solidly in
the Leibniz--l'H\^opital--Bernoulli--Euler school.  If so, the
question arises how modern commentators could assume that Cauchy meant
something else by \emph{infinitely small} than was customary in his
natural scientific milieu.

How was the notion perceived by Cauchy's contemporaries like Poisson
as well as a majority of Cauchy's colleagues at the \emph{\'Ecole
Polytechnique}?  Their comments (see below) indicate that their work
was a natural habitat for infinitesimals in the sense of the founders
of the calculus.  Thus, Cauchy's colleague Petit, a professor of
physics, requested that
\begin{quote}
this material [on differential calculus] be presented without certain
notions from algebra, which mainly had to do with series and which, he
alleged, the students would never have occasion to use in the
[engineering] services. Moreover, he insisted that the method of
\emph{infinitesimals} be used.  (Petit as translated by Belhoste in
\cite{Be91}, 1991, p.\;65; emphasis added)
\end{quote}
In a similar vein, de Prony reported:
\begin{quote}
I will finish my observations on the course in pure analysis by
manifesting the desire to see the use of the algorithm of imaginaries
[i.e., complex numbers] reduced to what is strictly necessary.  I have
been astonished, for instance, to see the expression of the element of
a curve, given in polar coordinates, derived [by Cauchy] from an
analysis using this algorithm; it follows much more quickly and with
greater ease from a consideration of \emph{infinitesimals}.  (de Prony
as translated in \cite[p.\;83]{Be91}; emphasis added)
\end{quote}
Poisson and de Prony both championed the use of infinitesimals through
their influence on the \emph{Conseil de Perfectionnement} (CP) of the
\emph{\'Ecole}, as noted by Gilain:
\begin{quote}
[O]n trouve dans le programme officiel, adopt\'e par le CP, une
modification significative : l'ajout dans les applications
g\'eom\'etriques du calcul diff\'erentiel et int\'egral, et dans le
programme de m\'ecanique, de l'instruction d'utiliser les
\emph{infiniment petits}.  M\^eme si l'auteur de cette proposition
n'est pas mentionn\'e dans les Proc\`es-verbaux, on peut penser
qu'elle \'emane des examinateurs de math\'e\-matiques, \emph{Poisson}
et \emph{de Prony}, qui animaient en g\'en\'eral la commission
programme du CP, et dont on conna\^\i t les convictions en faveur de
la m\'ethode des \emph{infiniment petits}.  (Gilain \cite{Gi89}, 1989,
\S 32; emphasis added)
\end{quote}
Like Cauchy, his contemporaries de Prony, Petit, and Poisson saw
infinitesimals as a natural tool both in teaching and in research,
though they were critical of what they saw as excessive rigor in
Cauchy's teaching.

Paying proper attention to the scientific context of the period goes
hand-in-hand with looking at Cauchy's practices and procedures on his
own terms, or as close as possible to his own terms, without
necessarily committing oneself to an Archimedean interpretation
thereof.  Thus, Ferraro writes:
\begin{quote}
Cauchy uses infinitesimal neighborhoods of~$x$ in a decisive way
{\ldots} Infinitesimals are not thought as a mere \emph{fa\c con de
parler}, but they are conceived as numbers, though a theory of
infinitesimal numbers is lacking.  (Ferraro \cite{Fe08}, 2008,
p.\;354)
\end{quote}
This comment by Ferraro is remarkable for two reasons:
\begin{enumerate}
\item
it displays a clear grasp of the distinction between procedure and
ontology (see B\l aszczyk et al.\;(\cite{17d}, 2017);
\item
it is a striking admission concerning the \emph{bona fide} nature of
Cauchy's infinitesimals.
\end{enumerate}
Ferraro's comment is influenced by Laugwitz's perceptive analysis of
Cauchy's sum theorem in (\cite{La87}, 1987), a paper cited several
times on Ferraro's page 354.

\section{Modern infinitesimals in relation to Cauchy's procedures}
\label{s71}

While set-theoretic justifications for a modern framework,
Archi\-medean or otherwise, are obviously not to be found in Cauchy,
Cauchy's \emph{procedures} exploiting infinitesimals find closer
proxies in Robinson's framework for analysis with infinitesimals than
in a Weierstrassian framework.  In this section we outline a
construction of a hyperreal extension
\begin{equation}
\label{e71}
\R\hookrightarrow\astr,
\end{equation}
and point out specific similarities between procedures using the
hyperreals, on the one hand, with Cauchy's procedures, on the other.

Let~$\R^{\N}$ denote the ring of sequences of real numbers, with
arithmetic operations defined termwise.  Then we have
\[
\astr=\R^{\N}\!/\,\text{MAX}
\]
where MAX is the maximal ideal consisting of all ``negligible''
sequences~$(u_n)$.  Here a sequence is negligible if it vanishes for a
set of indices of full measure~$\xi$, namely,~$\xi\big(\{n\in\N\colon
u_n=0\}\big)=1$.  Here
\[
\xi\colon \mathcal{P}(\N)\to \{0,1\}
\]
is a finitely additive probability measure taking the value~$1$ on
cofinite sets, where~$\mathcal{P}(\N)$ is the set of subsets of~$\N$.
The subset~$\mathcal{F}_\xi\subseteq\mathcal{P}(\N)$ consisting of
sets of full measure~$\xi$ is called a free ultrafilter.  These
originate with Tarski (\cite{Ta30}, 1930).  The set-theoretic
presentation of an infinitesimal-enriched continuum was therefore not
available prior to that date.  

The embedding $\eqref{e71}$ uses constant sequences.  We can therefore
define the subring
\begin{equation}
\label{s72}
\hr\subseteq\astr
\end{equation}
to be the set of the finite elements of~$\astr$; i.e., elements
smaller in absolute value than some real number, relying on the
embedding~\eqref{e71}.  The subring~\eqref{s72} admits a map~$\st$
to~$\R$, known as \emph{standard part}
\begin{equation}
\label{e61}
\st\colon \hr\to\R,
\end{equation}
which rounds off each finite hyperreal number to its nearest real
number.  This enables one, for instance, to define continuity and the
derivative as follows.  Following Robinson, we say that a function
\begin{quote}
$f(x)$ is continuous in [an open interval]~$(a,b)$ if 
\begin{equation}
\label{e74}
f(x_0+\eta)=_1 f(x_0)
\end{equation}
for all \emph{standard}~$x_0$ in the open interval and for all
infinitesimal~$\eta$.  (Robinson \cite{Ro61}, 1961, p.\;436; emphasis
in the original; counter~\eqref{e74} added)
\end{quote}
Robinson's symbol ``$=_1$'' denotes the relation of infinite
proximity.  Robinson's notation $=_1$ in \cite{Ro61} for infinite
proximity was replaced by $\simeq$ in his books and by $\approx$ in
most modern sources in infinitesimal analysis.

We also define the derivative of $t=f(s)$ as
\begin{equation}
\label{e62}
f'(s)=\st\left(\frac{\Delta t}{\Delta s}\right)
\end{equation}
where~$\Delta s\ne0$, or equivalently~$f'(s)$ is the standard real
number such that
\begin{equation}
\label{e75}
f'(s) \approx \frac{\Delta t}{\Delta s}.
\end{equation}
Such a definition parallels Cauchy's definition \eqref{e31b} of
derivative, more closely than any \emph{Epsilontik} definition.  Limit
is defined in terms of standard part, e.g., by setting
\begin{equation}
\label{e83}
\lim_{s\to0}f(s)=\st(f(\epsilon))
\end{equation}
where~$\epsilon$ is a nonzero infinitesimal.  This definition of limit
via the standard part is analogous to Cauchy's \emph{limite},
similarly defined in terms of infinitesimals, as analyzed in
Section~\ref{s21}.  For more details on Robinson's framework see e.g.,
Fletcher et al.\;(\cite{17f}, 2017).

\section{Conclusion}

We have argued that Cauchy's work on integral geometry, centers of
curvature, and other applications exploits infinitesimals as atomic
entities not reducible to simpler ones (such as terms in a sequence).

The oft-repeated claim, as documented e.g., in Bair et
al.\;(\cite{17a}, 2017) and Bascelli et al.\;(\cite{18e}, 2018), that
``Cauchy's infinitesimal is a variable with limit~$0$'' is a
reductionist view of Cauchy's foundational stance, at odds with much
compelling evidence in Cauchy's writings, as we argued in
Sections~\ref{s43} through \ref{six}.  Cauchy's notion of
infinitesimal was therefore close to that of his contemporary
scientists including Poisson, as we saw in Section~\ref{s64}.

While Cauchy did give an occasional \emph{Epsilontik} proof that today
would be interpreted in the context of an Archimedean continuum, his
techniques relying on infinitesimals find better proxies in a modern
framework exploiting an infinitesimal-enriched continuum.  Cauchy's
infinitesimal techniques in fields as diverse as geometric
probability, differential geometry, continuity and convergence are
just as viable as his \emph{Epsilontik} techniques.

Robinson first proposed an interpretation of Cauchy's
\emph{procedures} in the framework of a modern theory of
infinitesimals in (\cite{Ro66}, 1966, chapter\;10).  A
\emph{set-theoretic foundation} for infinitesimals could not have been
provided by Cauchy for obvious reasons, but Cauchy's \emph{procedures}
find closer proxies in modern infinitesimal frameworks than in modern
Archi\-medean ones.

\section*{Acknowledgments} 

We thank Peter Fletcher, El\'\i as Fuentes Guill\'en, Karel Hrbacek,
Taras Kudryk, David Pierce, and David Sherry for helpful suggestions.
V.\;Kanovei was partially supported by RFBR Grant 17-01-00705.

\end{document}